\documentclass[12pt]{article}

\usepackage{amsthm,amsfonts,amsbsy,amssymb,amsmath,wrapfig}

\newcommand{\ba}[1]{\bigl\langle#1\bigr\rangle}

\DeclareMathOperator{\Supp}{Supp} \DeclareMathOperator{\Ker}{Ker}

\newtheorem{thm}{Theorem}

\newtheorem{prop}{Proposition}

\theoremstyle{definition}

\newtheorem{dfn}{Definition}
\newtheorem*{dfn*}{Definition}

\theoremstyle{remark}

\newtheorem{rem}{Remark}

\DeclareMathOperator{\otimesZ}{\otimes\hspace{1pt}\rule[-3pt]{0pt}{0pt}_{\mathbb{Z}}}

\begin{document}

\begin{center}
{\large \bfseries On solvable spherical subgroups

of semisimple algebraic groups}

\vskip10pt

{\scshape Roman Avdeev}

\vskip10pt

\end{center}

This report contains a new approach to classification of solvable
spherical subgroups of semisimple algebraic groups. This approach is
completely different from the approach by D. Luna \cite{Lu} and
provides an explicit classification.

Let $G$ be a connected semisimple complex algebraic group. We fix a
Borel subgroup $B \subset G$ and a maximal torus $T \subset B$. We
denote by $U$ the maximal unipotent subgroup of $G$ contained in
$B$. The Lie algebras of $G$, $B$, $U$, $\ldots$ are denoted by
$\mathfrak g$, $\mathfrak b$, $\mathfrak u$, $\ldots$, respectively.
Let $\Delta \subset \mathfrak X(T)$ be the root system of $G$ with
respect to $T$ (where $\mathfrak X(T)$ denotes the character lattice
of $T$). The subsets of positive roots and simple roots with respect
to $B$ are denoted by $\Delta_+$ and $\Pi$, respectively. For any
root $\alpha \in \Delta_+$ consider its expression of the form
$\alpha = \sum\limits_{\gamma \in \Pi} k_\gamma \gamma$. We put
$\Supp \alpha = \{\gamma \in \Pi \mid k_\gamma > 0\}$. The root
subspace of the Lie algebra $\mathfrak g$ corresponding to a root
$\alpha$ is denoted by $\mathfrak g_\alpha$. The symbol $\ba{A}$
will denote the linear span of a subset $A \subset \mathfrak X(T)$
in $\mathfrak X(T) \otimesZ \mathbb Q$.

Let $H \subset B$ be a connected solvable algebraic subgroup of $G$
and $N \subset U$ its unipotent radical. We say that the group $H$
is \emph{standardly embedded} into $B$ (with respect to $T$) if the
subgroup $S = H \cap T$ is a maximal torus of $H$. Obviously, in
this case we have $H = S \rightthreetimes N$.

Suppose that $H \subset G$ is a connected solvable subgroup
standardly embedded into $B$. Then we can consider the natural
restriction map $\tau: \mathfrak X(T) \to \mathfrak X(S)$. Put $\Phi
= \tau(\Delta_+) \subset \mathfrak X(S)$. Then we have $\mathfrak u
= \bigoplus \limits_{\lambda \in \Phi} \mathfrak u_\lambda$, where
$\mathfrak u_\lambda \subset \mathfrak u$ is the weight subspace of
weight $\lambda$ with respect to $S$. Similarly, we have $\mathfrak
n = \bigoplus \limits_{\lambda \in \Phi} \mathfrak n_\lambda$, where
$\mathfrak n_\lambda = \mathfrak u_\lambda \cap \mathfrak n \subset
\mathfrak u_\lambda$. Denote by $c_\lambda$ the codimension of
$\mathfrak n_\lambda$ in~$\mathfrak u_\lambda$.

Recall that a subgroup $H\subset G$ is called \emph{spherical} if
the group $B$ has an open orbit in $G/H$. The following theorem
provides a convenient criterion of sphericity for connected solvable
subgroups of $G$.

\begin{thm}
Suppose $H \subset G$ is a connected solvable subgroup standardly
embedded into $B$. Then the following conditions are equivalent:

\textup{(1)} $H$ is spherical in $G$;

\textup{(2)} $c_\lambda \leqslant 1$ for every $\lambda \in \Phi$,
and the weights $\lambda$ with $c_\lambda = 1$ are linearly
independent in~$\mathfrak X(S)$.
\end{thm}

Now we suppose that $H$ is spherical.

\begin{dfn}
A root $\alpha \in \Delta_+$ is called \emph{marked} if $\mathfrak
g_\alpha \not\subset \mathfrak n$.
\end{dfn}

An important role of marked roots in studying solvable spherical
subgroups is clear from the theorem below.

\begin{thm}
Up to conjugation by elements of\: $T$, the subgroup $H$ is uniquely
determined by its maximal torus $S \subset T$ and the set $\Psi
\subset \Delta_+$ of marked roots.
\end{thm}

\begin{rem}
The subgroup $H$ is explicitly recovered from $S$ and $\Psi$.
\end{rem}

Marked roots have the following property: if $\alpha \in \Psi$ and
$\alpha = \beta + \gamma$ for some roots $\beta, \gamma \in
\Delta_+$, then exactly one of two roots $\beta, \gamma$ is marked.
Taking this property into account, we say that a marked root $\beta$
is \emph{subordinate} to a marked root $\alpha$, if $\alpha = \beta
+ \gamma$ for some root $\gamma \in \Delta_+$. Given a marked root
$\alpha$, we denote by $C(\alpha)$ the set consisting of $\alpha$
and all marked roots subordinate to $\alpha$. Further, we say that a
marked root $\alpha$ is \emph{maximal} if it is not subordinate to
any other marked root. Let $\mathrm M$ denote the set of all maximal
marked roots.

\begin{prop}\label{associated_simple_root}
For any marked root $\alpha$ there exists a unique simple root
$\pi(\alpha) \in \Supp \alpha$ with the following property: if
$\alpha = \beta + \gamma$ for some roots $\beta, \gamma \in
\Delta_+$, then the root $\beta$ is marked iff $\pi(\alpha) \notin
\Supp \beta$ (and so the root $\gamma$ is marked iff $\pi(\alpha)
\notin \Supp \gamma$).
\end{prop}

\begin{dfn}
If $\alpha$ is a marked root, then the simple root $\pi(\alpha)$
appearing in Proposition~\ref{associated_simple_root} is called the
\emph{simple root associated with the marked root} $\alpha$.
\end{dfn}

From Proposition~\ref{associated_simple_root} we see that for any
marked root $\alpha$ the set $C(\alpha)$ is uniquely determined by
the simple root $\pi(\alpha)$. Therefore, the whole set $\Psi$ is
uniquely determined by the set $\mathrm M$ and the map $\pi\colon
\mathrm M \to \Pi$.

\begin{thm}[Classification of marked roots]
All possibilities for a marked root $\alpha$ and the simple root
$\pi(\alpha)$ are presented in Table~\ref{table_marked_roots}.
\end{thm}

\textbf{Notations used in Table~\ref{table_marked_roots}}. The
symbol $\Delta(\alpha)$ denotes the root system generated by $\Supp
\alpha$, i.~e. $\Delta(\alpha) = \ba{\Supp \alpha} \cap \Delta$. We
suppose that $\Supp \alpha = \{\alpha_1, \ldots, \alpha_n\}$. The
numeration of simple roots of simple Lie algebras is the same as
in~\cite{VO}.

\begin{table}[h]

\caption{Marked roots}\label{table_marked_roots}

\begin{center}

\begin{tabular}{|c|c|c|c|}

\hline

¹ & type of $\Delta(\alpha)$ & $\alpha$ & $\pi(\alpha)$ \\

\hline

1 & any of rank $n$ & $\alpha_1 + \alpha_2 + \ldots + \alpha_n$ &
$\alpha_1, \alpha_2, \ldots, \alpha_n$\\

\hline

2 & $\mathrm B_n$ & $\alpha_1 + \alpha_2 + \ldots + \alpha_{n-1} +
2\alpha_n$ &
$\alpha_1, \alpha_2, \ldots, \alpha_{n-1}$\\

\hline

3 & $\mathrm C_n$ & $2\alpha_1 + 2\alpha_2 + \ldots + 2\alpha_{n-1}
+ \alpha_n$
& $\alpha_n$\\

\hline

4 & $\mathrm F_4$ & $2\alpha_1 + 2\alpha_2 + \alpha_3 + \alpha_4$ &
$\alpha_3, \alpha_4$\\

\hline

5 & $\mathrm G_2$ & $2\alpha_1 + \alpha_2$ & $\alpha_2$\\

\hline

6 & $\mathrm G_2$ & $3\alpha_1 + \alpha_2$ & $\alpha_2$\\

\hline
\end{tabular}

\end{center}

\end{table}

For any root $\alpha \in \Delta_+$ consider the (connected) Dynkin
diagram $D(\alpha)$ of the set $\Supp \alpha$.  We say that a root
$\delta \in \Supp \alpha$ is \emph{terminal with respect to $\Supp
\alpha$} if the node of $D(\alpha)$ corresponding to $\delta$ is
connected by an edge (possibly, multiple) with exactly one other
node of $D(\alpha)$.

Now we introduce some conditions on a pair of marked roots $(\alpha,
\beta)$. These
\begin{wrapfigure}{r}{150\unitlength}

\begin{picture}(150,115)
\put(70,75){\circle{4}} \put(65,83){$\gamma_0$}

\put(70,73){\line(0,-1){18}} \put(70,53){\circle{4}}
\put(76,51){$\gamma_1$} \put(70,51){\line(0,-1){11}}
\put(70,36){\circle*{0.5}} \put(70,32){\circle*{0.5}}
\put(70,28){\circle*{0.5}} \put(70,25){\line(0,-1){11}}
\put(70,12){\circle{4}} \put(76,10){$\gamma_r$}

\put(72,76){\line(2,1){16}} \put(68,76){\line(-2,1){16}}

\put(90,85){\circle{4}} \put(88,74){$\beta_1$}
\put(92,86){\line(2,1){10}} \put(106,93){\circle*{0.5}}
\put(110,95){\circle*{0.5}} \put(114,97){\circle*{0.5}}
\put(118,99){\line(2,1){10}} \put(130,105){\circle{4}}
\put(128,94){$\beta_q$}

\put(50,85){\circle{4}} \put(45,74){$\alpha_1$}
\put(48,86){\line(-2,1){10}} \put(34,93){\circle*{0.5}}
\put(30,95){\circle*{0.5}} \put(26,97){\circle*{0.5}}
\put(22,99){\line(-2,1){10}} \put(10,105){\circle{4}}
\put(05,94){$\alpha_p$}

\end{picture}

{

\makeatletter

\renewcommand{\abovecaptionskip}{0pt}
\renewcommand{\belowcaptionskip}{0pt}

\renewcommand{\@makecaption}[2]{
\vspace{\abovecaptionskip}%
\sbox{\@tempboxa}{#1. #2}%
\global\@minipagefalse \hbox to \hsize {{\hfil #1. #2\hfil}}
\vspace{\belowcaptionskip}}

\makeatother

\caption{}\label{diagram_difficult}

}

\end{wrapfigure} conditions will be used later.

$(\mathrm D0)$ $\Supp \alpha \cap \Supp \beta = \varnothing$

$(\mathrm D1)$ $\Supp \alpha \cap \Supp \beta = \{\delta\}$,
$\delta$ is terminal with respect to both $\Supp \alpha$ and $\Supp
\beta$, $\pi(\alpha) \ne \delta \ne \pi(\beta)$

$(\mathrm D2)$ the Dynkin diagram of the set $\Supp \alpha \cup
\Supp \beta$ has the form shown on Figure~\ref{diagram_difficult}
(for some $p, q, r \geqslant 1$), $\alpha = \alpha_1 + \ldots +
\alpha_p + \gamma_0 + \gamma_1 + \ldots + \gamma_r$, $\beta =
\beta_1 + \ldots + \beta_q + \gamma_0 + \gamma_1 + \ldots +
\gamma_r$, $\pi(\alpha) \notin \Supp \alpha \cap \Supp \beta$,
$\pi(\beta) \notin \Supp \alpha \cap \Supp \beta$

$(\mathrm E1)$ $\Supp \alpha \cap \Supp \beta = \{\delta\}$,
$\delta$ is terminal with respect to both $\Supp \alpha$ and $\Supp
\beta$, $\delta = \pi(\alpha) = \pi(\beta)$, $\alpha - \delta \in
\Delta_+$, $\beta - \delta \in \Delta_+$

$(\mathrm E2)$ the Dynkin diagram of the set $\Supp \alpha \cup
\Supp \beta$ has the form shown on Figure~\ref{diagram_difficult}
(for some $p, q, r \geqslant 1$), $\alpha = \alpha_1 + \ldots +
\alpha_p + \gamma_0 + \gamma_1 + \ldots + \gamma_r$, $\beta =
\beta_1 + \ldots + \beta_q + \gamma_0 + \gamma_1 + \ldots +
\gamma_r$, $\pi(\alpha) = \pi(\beta) \in \Supp \alpha \cap \Supp
\beta$

Next, we need to introduce an equivalence relation on $\mathrm M$.
For any two roots $\alpha, \beta \in \mathrm M$ we write $\alpha
\sim \beta$ iff $\tau(\alpha) = \tau(\beta)$. Having introduced this
equivalence relation, to each connected solvable spherical subgroup
$H$ standardly embedded into $B$ we assign the set of combinatorial
data $(S, \mathrm M, \pi, \sim)$.

\begin{thm}\label{uniqueness_theorem}
The above assignment is a one-to-one correspondence between the
following two sets:

\textup{(1)} the set of all connected solvable spherical subgroups
standardly embedded into $B$, up to conjugation by elements of $T$;

\textup{(2)} the set of all sets $(S, \mathrm M, \pi, \sim)$, where
$S \subset T$ is a torus, $\mathrm M$ is a subset of $\Delta_+$,
$\pi\colon \mathrm M \to \Pi$ is a map, $\sim$ is an equivalence
relation on $\mathrm M$, and the set $(S, \mathrm M, \pi, \sim)$
satisfies the following conditions:

$(\mathrm A)$ $\pi(\alpha) \in \Supp \alpha$ for any $\alpha \in
\mathrm M$, and the pair $(\alpha, \pi(\alpha))$ is contained in
Table~\ref{table_marked_roots};

$(\mathrm D)$ if $\alpha, \beta \in \mathrm M$ and $\alpha \not\sim
\beta$, then one of conditions $(\mathrm D0)$, $(\mathrm D1)$,
$(\mathrm D2)$ holds;

$(\mathrm E)$ if $\alpha, \beta \in \mathrm M$ and $\alpha \sim
\beta$, then one of conditions $(\mathrm D0)$, $(\mathrm D1)$,
$(\mathrm E1)$, $(\mathrm D2)$, $(\mathrm E2)$ holds;

$(\mathrm C)$ for any $\alpha \in \mathrm M$ the condition $\Supp
\alpha \not\subset \bigcup\limits_{\beta \in \mathrm M \backslash
\{\alpha\}}\Supp \beta$ holds;

$(\mathrm T)$ $\left.\Ker \tau\right|_R = \ba{\alpha - \beta \mid
\alpha, \beta \in \mathrm M, \alpha \sim \beta}$, where $R =
\ba{\bigcup\limits_{\gamma \in \mathrm M} \Supp \gamma}$.
\end{thm}

\begin{rem}
The unipotent radical $N\subset U$ of a connected solvable spherical
subgroup $H$ standardly embedded into $B$ is uniquely (up to
conjugation by elements of $T$) determined by the set $(\mathrm M,
\pi, \sim)$ satisfying conditions $(\mathrm A)$, $(\mathrm D)$,
$(\mathrm E)$, $(\mathrm C)$.
\end{rem}

To complete the classification of connected solvable spherical
subgroups of $G$ up to conjugation, it remains to determine all sets
of combinatorial data that correspond to one conjugacy class of such
subgroups. Consider again a connected solvable spherical subgroup $H
\subset G$ standardly embedded into $B$. We say that a marked root
$\alpha$ is \emph{regular} if the projection of $\mathfrak n$ to the
root space $\mathfrak g_\alpha$ is zero. Choose any regular marked
simple root $\alpha \in \Delta_+$ and fix an element $n_\alpha \in
N_G(T)$ such that its image in the Weyl group $W$ is the simple
reflection $r_\alpha$ corresponding to $\alpha$ (here $N_G(T)$ is
the normalizer of $T$ in $G$). Obviously, the group $n_\alpha H
n_\alpha^{-1}$ is also standardly embedded into $B$. Its set of
combinatorial data is $(n_\alpha S n_\alpha^{-1}, \mathrm M', \pi',
\sim')$ for some $\mathrm M'$, $\pi'$, and $\sim'$. In order to
determine $\mathrm M'$, $\pi'$, and $\sim'$, we consider two cases:

(1) if $\alpha \in \Supp \delta$ for some $\delta \in
r_\alpha(\mathrm M \backslash \{\alpha\})$, then $\mathrm M' =
r_\alpha(\mathrm M \backslash \{\alpha\})$, $\pi'(\beta) =
\pi(r_\alpha(\beta))$ for any $\beta \in \mathrm M'$, $\beta \sim'
\gamma$ iff $r_\alpha(\beta) \sim r_\alpha(\gamma)$ for any $\beta,
\gamma \in \mathrm M'$;

(2) if $\alpha \notin \Supp \delta$ for any $\delta \in
r_\alpha(\mathrm M \backslash \{\alpha\})$, then $\mathrm M' =
r_\alpha(\mathrm M \backslash \{\alpha\}) \cup \{\alpha\}$,
$\pi'(r_\alpha(\beta)) = \pi(\beta)$ for any $\beta \in \mathrm
M\backslash\{\alpha\}$, $\pi'(\alpha) = \alpha$, $r_\alpha(\beta)
\sim' r_\alpha(\gamma)$ iff $\beta \sim \gamma$ for any $\beta,
\gamma \in \mathrm M \backslash \{\alpha\}$, $\alpha \not\sim'
\beta$ for any $\beta \in \mathrm M' \backslash\{\alpha\}$.

Transformations of the form $H \mapsto n_\alpha H n_\alpha^{-1}$
described above are called \emph{elementary transformations}.

\begin{thm} \label{conj}
Suppose $H_1, H_2 \subset G$ are two connected solvable spherical
subgroups standardly embedded into $B$, and $H_2 = g H_1 g^{-1}$ for
some $g \in G$. Then:

\textup{(1)} $H_2 = \sigma H_1 \sigma^{-1}$ for some $\sigma \in
N_G(T)$;

\textup{(2)} $H_2$ can be obtained from $H_1$ by applying a finite
sequence of elementary transformations.
\end{thm}

Thus, Theorems~\ref{uniqueness_theorem} and~\ref{conj} give a
complete classification of connected solvable spherical subgroups of
$G$.

\end{document}